\documentclass[11pt,a4paper]{article}
\usepackage{amssymb,amsmath,amsthm,amscd,graphicx}
\usepackage{amsmath,amsfonts,amssymb}
\newtheorem{theorem}{Theorem}[section]
\newtheorem{proposition}{Propostion}[section]

\newtheorem{lemma}{Lemma}[section]
\newtheorem{definition}{Definition}[section]
\newtheorem{remark}{Remark}[section]

\begin{document}
\title{Lattice Approximations of Reflected Stochastic Partial Differential Equations Driven by Space-Time  White Noise}
\author{\
Tusheng Zhang$^{1}$}

\footnotetext[1]{\ School of Mathematics, University of Manchester,
Oxford Road, Manchester M13 9PL, England, U.K. Email:
tusheng.zhang@manchester.ac.uk}
 \maketitle{} {\bf
Abstract}: We introduce a discretization/approximation scheme for reflected
stochastic partial differential equations driven by space-time white noise through
systems of reflecting stochastic differential equations.
To establish the convergence of the scheme, we study the existence and uniqueness of solutions of Skorohod-type deterministic systems
on time-dependent domains. We also need to establish the convergence of an approximation scheme for deterministic parabolic obstacle problems. Both are of independent interest on their own.


\section{Introduction}
\setcounter{equation}{0} Consider the following stochastic partial
differential equation (SPDE) with reflection:
\begin{eqnarray}\label{original equation}\left\{
\begin{array}{l}\frac{\partial u(t,x)}{\partial t}-\frac{\partial^2 u(t,x)}{\partial x^2}=f(t,x,u(t,x))+\sigma(t,x,u(t,x))\dot{W}(t,x)+\eta;\\
u(0,\cdot)=u_0, \quad u(t,x)\geq 0;\\
 u(t,0)=u(t,1)=0,
\end{array}
\right.
\end{eqnarray}
where $\dot{W}$ denotes the space-time white noise defined on a
complete probability space $(\Omega,
\mathcal{F},\{\mathcal{F}_t\}_{t\geq 0 },P)$, where
$\mathcal{F}_t=\sigma(W(s,x): x\in [0,1], 0\leq s\leq t)$; $u_0$ is
a non-negative  continuous function on $[0,1]$, which vanishes at 0
and 1;  $\eta(t,x)$ is a random measure which is a part of the
solution pair $(u, \eta)$ and plays the role of a local time that prevents
 the solution $u$ from being negative. The coefficients $f$ and $\sigma$ are
measurable mappings from $\mathbb{R}_+\times [0,1]\times
\mathbb{R}$ into $\mathbb{R}$. Let $C^2_0([0,1])$ denote the space of twice differentiable functions $\phi$
on $[0, 1]$ satisfying $\phi(0)=\phi(1)=0$.  The following
definition is taken from \cite{MP}, \cite{NP}.
\begin{definition} A pair $(u,\eta)$ is said to
be a solution of equation (\ref{original equation}) if\\
(i) $u$ is a continuous random field on $\mathbb{R}_+\times [0,1];$
$u(t,x)$ is $\mathcal{F}_t$ measurable and $u(t,x)\geq 0\ a.s.$\\
(ii) $\eta$ is a random measure on $\mathbb{R}_+\times (0,1)$ such
that\\
\indent (a) $\eta(\{t\}\times (0,1))=0,\quad \forall t\geq 0.$\\
\indent (b) $\int_0^t\int_0^1 x(1-x)\eta(ds,dx)<\infty, t\geq 0.$\\
\indent (c) $\eta$ is adapted in the sense that for any measurable
mapping $\psi:$
$$\int_0^t\int_0^1\psi(s,x)\eta(ds,dx)\ \mbox{is}\ \mathcal{F}_t\ \mbox{measurable}.$$
(iii) $\{u, \eta\}$ solves the parabolic SPDE in the following sense (
$(\cdot,\cdot)$ denotes the scalar product in $L^2[0,1]$): $\forall
t\in \mathbb{R}_+, \phi\in C^2_0([0,1])$ with $\phi(0)=\phi(1)=0,$
\begin{eqnarray*}
&&(u(t),\phi)-\int_0^t(u(s),\phi'')ds-\int_0^t(f(s,\cdot,u(s)),\phi)ds\\
&=&(u_0,
\phi)+\int_0^t\int_0^1\phi(x)\sigma(s,x,u(s,x))W(ds,dx)+\int_0^t\int_0^1\phi(x)\eta(ds,dx)\
a.s,\end{eqnarray*} where $u(t):=u(t,\cdot)$.

\noindent(iv) $\int_Q ud\eta=0,$\quad where $Q=
\mathbb{R}_+\times (0,1).$
\end{definition}
The SPDEs with reflection driven by space-time white noise was first studied by Nualart and Pardoux in \cite{NP}
(PTRF 1992) when $\sigma(\cdot)=1$, and by Donati-Martin and Pardoux
 in \cite{MP}(in PTRF 1993) for general diffusion coefficient $\sigma$. The uniqueness of the solution and large deviations
 were obtained by Xu and Zhang in \cite{XZ}.
 SPDEs with reflection can  be used to
model the evolution of random interfaces near a hard wall. It was
proved by T. Funaki and S. Olla in \cite{FO} that the fluctuations
of a $\nabla \phi$ interface model near a hard wall converge in law to the stationary solution of a SPDE with reflection.
Various properties of the solution of equation \eqref{original
equation} were studied since then. The hitting properties were investigated by Dalang, Mueller and Zambotti in \cite{DMZ}.
Integration by parts formulae associated with SPDEs with reflection, occupation densities were established by Zambotti in \cite{ZL-1},\cite{ZL-2} and \cite{ZL-3}. The strong Feller properties and the large deviations for invariant measures were studied by Zhang in \cite{Z-1}, \cite{Z-2}.

\vskip 0.4cm

 \indent The purpose of this paper is to develop a numerical scheme(particle system approximations) for the reflected
 stochastic partial differential equations. This is a challenging problem which has been open for some time.
 Part of the difficulties is caused by the singularities of the space-time white noise. For example, Ito formula is not available for
  this  type of equations. Part of the difficulties  lie in the discretization of the random measure
 $\eta$ appeared in the equation (\ref{original equation}). We introduce a discretization scheme through systems of reflecting stochastic differential equations. As the dimensions of the reflecting systems tend to infinity, the problem is to compare and control the systems with different dimensions. To this end  we study Skorohod-tye deterministic problems on time-dependent domains and prove a useful a priori estimate for the solutions in terms  of time-dependent boundaries. To prove the convergence of the scheme, we also need to establish the convergence of a discretization scheme of deterministic parabolic obstacle problems. These preliminary results are of independent interest.

\vskip 0.3cm
The discretization scheme for stochastic heat  equations driven by space-time white noise was first introduced by Gy\"o{}ngy in \cite{G}, \cite{G-1}. Approximation scheme for SPDEs of elliptic type was discussed by Martinez and Sanz-Sol\`e in \cite{MS}. Discretisations for stochastic wave equations were investigated by Quer-Sardanyons and Sanz-Sol\`e in \cite{QS}. Numerical schemes for stochastic evolution equations were obtained by Gy\"o{}ngy and Millet in \cite{GM}.

\vskip 0.4cm
 \indent Let us now describe the content of the paper in more
 detail. In Section 2, we introduce the discretization scheme and the main result.
 Section 3 is to  study Skorohod -type problems on time-dependent domains in Euclidean spaces. We establish
 the existence and uniqueness  of the solution of the Skorohod type problem on domains with boundaries being
 continuous functions of time. We provide a bound of the solution in terms of the  boundaries of the domains, which plays an important role in the rest of the paper. In Section 4,  we introduce a discretization scheme for deterministic parabolic obstacle problems. We establish the convergence of the scheme first for smooth obstacles. In this case, we are able to show that the measure $\eta(dt,dx)$ appeared in the obstacle problem is absolutely continuous with respect to the Lebesgue measure $dt\times dx$ and the tightness of the approximating solutions. We prove the convergence of the scheme by identifying any limit of the approximating solutions as the unique solution of the parabolic obstacle problem.  We then extend the scheme  for continuous obstacles using the  a priori estimate obtained in Section 3 for Skorohod-type problems. The Section 5 is devoted to the proof of the convergence of the discretization scheme for SPDEs with reflection. We first relate the SPDEs with reflection to a random parabolic obstacle problem. We obtain the convergence of the scheme by carefully comparing it with the discretization scheme introduced for obstacle problems in Section 4. Here the results in Section 3  and Garsia Lemma for random fields will play an important role.
 \vskip 0.4cm
 \begin{remark}
 In this paper we discretize only the space variable using systems of reflecting stochastic differential equations (SDEs). Now there is known procedure to further discretize the reflected SDEs (see \cite{S}). Combining these two one could get the discretization for SPDEs with reflection both in time and space directions.
 \end{remark}

\section {The discretization scheme and the main result}
\setcounter{equation}{0}
\noindent We first introduce the conditions on the
coefficients. Let $f,\sigma$ are two measurable mappings
$$f,\sigma: \mathbb{R}_+\times [0,1]\times \mathbb{R}\rightarrow \mathbb{R}$$
satisfying:\\
{\bf (H.1)}. For any $T>0$, there exists a constant
$c(T)$  such that for any $x,y\in [0,1], t\in
[0,T]$, $u,v\in \mathbb{R}$,
\begin{eqnarray}\label{locally lipschitz}
|f(t,x,u)-f(t,y,v)|+|\sigma(t,x,u)-\sigma(t,y,v)|\leq
c(T)[|x-y|+|u-v|].
\end{eqnarray}
{\bf (H.2)}. For any $T>0$, there exists a constant
$c(T)$  such that for any $x\in [0,1], t\in
[0,T]$, $u\in \mathbb{R}$,
\begin{eqnarray}\label{lineargrowth}
|f(t,x,u)|\leq c(T)(1+|u|).
\end{eqnarray}
\vskip 0.4cm
For every integer $n\geq 1$ and $x=\frac{k}{n}, k=1,2,...,n-1$, define the processes $u^n(t, \frac{k}{n}), k=1,...,n-1$ as the solution
of the system of reflecting stochastic differential equations
\begin{eqnarray}\label{2.1}
du^n(t, \frac{k}{n})&=&n^2\left ( u^n(t, \frac{k+1}{n})-2u^n(t, \frac{k}{n})+u^n(t, \frac{k-1}{n})\right )dt\nonumber\\
&+&n\sigma (t, \frac{k}{n}, u^n(t, \frac{k}{n}))d\left (W(t,\frac{k+1}{n})-W(t,\frac{k}{n})\right )\nonumber\\
&+&f(t, \frac{k}{n}, u^n(t, \frac{k}{n}))dt+d\eta_k^{n}(t)\\
u^n(t, \frac{k}{n})&\geq &0, \quad k=1,2,...,n-1,\nonumber\\
u^n(t,0)&=&u^n(t,1)=0, \quad k=1,2,...,n-1,\nonumber
\end{eqnarray}
with initial condition
\begin{equation}\label{2.2}
u^n(0,\frac{k}{n})=u_0(\frac{k}{n}), \quad  k=1,...,n-1.
\end{equation}
\begin{definition}
We say that $\{ u^n(t, \frac{k}{n}), \eta_k^{n}(t), k=1,...,n-1\}$ is a solution to the reflecting system (\ref{2.1}) if

(i) for every $k\geq 1$, $u^n(t, \frac{k}{n}), t\geq 0$ is an adapted,  non-negative, continuous process,

(ii) for $k\geq 1$,  $\eta_k^{n}(t), t\geq 0$ is an adapted, continuous, increasing process with $\eta_k^{n}(0)=0$,

(iii) for every $t\geq 0$ and $k=1,2,...,n-1$,
\begin{eqnarray}\label{2.3}
u^n(t, \frac{k}{n})&=&u_0(\frac{k}{n})+n^2\int_0^t \left ( u^n(s, \frac{k+1}{n})-2u^n(s, \frac{k}{n})+u^n(s, \frac{k-1}{n})\right )ds\nonumber\\
&+&\int_0^tf(s, \frac{k}{n}, u^n(s, \frac{k}{n}))ds+\eta_k^{n}(t) \nonumber\\
&+&\int_0^tn\sigma (s, \frac{k}{n}, u^n(s, \frac{k}{n}))d\left (W(s,\frac{k+1}{n})-W(s,\frac{k}{n})\right )
\end{eqnarray}
almost surely,

(iv) $\int_0^t u^n(s, \frac{k}{n})\eta_k^{n}(ds)=0$, for all $t\geq 0$, \quad\quad  $k=1,2,...,n-1$.
\end{definition}
\vskip 0.4cm

Set
$$u^n_k(t):=u^n(t, \frac{k}{n})$$
$$W^n_k(t)=\sqrt{n}\left ( W(t, \frac{k+1}{n})-W(t, \frac{k+1}{n})\right )$$
for $k=1,...,n-1$.
Let $A^n=(A_{ki}^n)$ denote the $(n-1)\otimes (n-1)$ matrix with elements $A_{kk}^n=-2$, $A_{ki}^n=1$ for $|k-i|=1$, $A_{ki}^n=0$ for $|k-i|>1$.
The system (\ref{2.1}) is regarded as a $(n-1)$-dimensional reflected SDE on the domain $D_n=\{(z_1,...,z_{n-1}); z_k\geq 0, k=1,...,n-1\}$
written as
\begin{eqnarray}\label{2.4}
du^n_k(t)&=&n^2\sum_{i=1}^{n-1} A^n_{ki}u^n_i(t)dt+f(t, \frac{k}{n}, u^n_k(t))dt \nonumber\\
&+&\sqrt{n}\sigma (t, \frac{k}{n}, u^n_k(t))dW^n_k(t)+d\eta_k^{n}(t)\\
u^n_k(0)&=&u_0(\frac{k}{n}), \quad k=1,2,...,n-1.\nonumber
\end{eqnarray}
As the domain $D_n$ is convex, the existence and uniqueness of the solution of the  system (\ref{2.4}) is  well known (see e.g. \cite{LS}).
\vskip 0.3cm
For every integer $n\geq 1$, define the random field
\begin{equation}\label{2.6}
u^n(t,x):=u^n(t, \frac{k}{n})+(nx-k)\left ( u^n(t, \frac{k+1}{n})-u^n(t, \frac{k}{n})\right )
\end{equation}
for $x\in [\frac{k}{n}, \frac{k+1}{n})$, $k=0,...,n-1$ with $u^n(t,0):=0$.
\vskip 0.3cm
The main result of the paper reads as

\begin{theorem}
Suppose (H.1) and (H.2) hold. Then for any $p\geq 1$, we have
\begin{equation}\label{2.7}
\lim_{n\rightarrow \infty}E[\sup_{0\leq t\leq T, 0\leq x\leq 1}|u^n(t,x)-u(t,x)|^p]=0.
\end{equation}
\end{theorem}
\vskip 0.4cm
We end this section with a description of the group generated by the matrix $A^n$.
For $j\geq 1$,  define the vector:
$$e_j=(e_j(k))=\left (\sqrt{\frac{2}{n}} \mbox{sin}(j\frac{1}{n}\pi),...,\sqrt{\frac{2}{n}} \mbox{sin}(j\frac{n-1}{n}\pi) \right ).$$
One can easily check that $\{e_j, j=1,...,n-1\}$ forms an orthonormal basis of $R^{n-1}$. Moreover, $e_j, j=1,...,n-1$ are eigenvectors of $n^2A^n$ with eigenvalues
$$\lambda^n_j:=-j^2\pi^2c^n_j,$$
where
$$\frac{4}{\pi^2}\leq c_j^n:=\frac{\mbox{sin}^2(\frac{j\pi}{2n})}{(\frac{j\pi}{2n})^2}\leq 1,$$
$j=1,...,n-1$.
Thus the group $G^n(t):=exp(n^2A^nt)$ generated by $n^2A^n$ on $R^{n-1}$ admits the following representation
\begin{equation}\label{2.5}
G^n(t)e=\sum_{k=1}^{n-1} e^{\lambda^n_kt}<e,e_k>e_k, \quad \quad e\in R^{n-1}.
\end{equation}

\section{ Deterministic Skorohod-type systems}
\setcounter{equation}{0}
In this section we study Skorohod-type problems on time dependent domains and obtain some a priori estimates.
\vskip 0.3cm
Set $a^+=a\vee 0$ and $a^-=(-a)\vee 0$ for $a\in R$. For a vector $b=(b_1,...,b_{n-1})\in R^{n-1}$, we will use the following notation
$$b^+=(b_1^+,...,b_{n-1}^+), \quad\quad b^-=(b_1^-,...,b_{n-1}^-).$$
It is clear that $b=b^+-b^-.$
\vskip 0.3cm
Note that $A^n$ introduced in Section 2 is a negative definite matrix. Furthermore, we also have the following
\begin{lemma}
It holds that
\begin{equation}\label{2.7}
<b^+, A^nb>\leq 0 \quad\quad \mbox{for all} \quad \quad b\in R^{n-1}.
\end{equation}
\end{lemma}
\vskip 0.4cm
\noindent{\bf Proof}. Write
\begin{equation}\label{2.8}
<b^+, A^nb>=<b^+, A^nb^+>-<b^+, A^nb^->.
\end{equation}
The first term on the right $<b^+, A^nb^+>$ is non-positive. Since $A^n_{ij}\geq 0$ for $i\not =j$,  we have
\begin{eqnarray*}
<b^+, A^nb^->&=&\sum_{i,j=1}^{n-1}b_i^+A_{ij}^nb_j^-\nonumber\\
&=&\sum_{i=1}^{n-1}b_i^+A_{ii}^nb_i^-+\sum_{i\not =j}b_i^+A_{ij}^nb_j^-\nonumber\\
&=&\sum_{i\not =j}b_i^+A_{ij}^nb_j^-\geq 0,
\end{eqnarray*}
(\ref{2.7}) follows. $\square$

\vskip 0.3cm
For $a=(a_1,...,a_{n-1}), b=(b_1,...,b_{n-1})\in R^{n-1}$, we write $a\geq b$ if $a_i\geq b_i$ for all $i=1,...,n-1$. Given $V=(V_1,...,V_{n-1})\in C([0,\infty )\rightarrow R^{n-1})$.
Consider  the following Skorohod-type problem with reflection in $R^{n-1}$:
\begin{eqnarray}\label{2.9}\left\{
\begin{array}{l}
dZ(t)=n^2A^nZ(t)dt+d\eta(t);\\
Z(t)\geq -V(t);\\
\int_0^T<Z(t)+V(t),d\eta(t)>=0, \quad\quad \mbox{for all} \quad T>0.
\end{array}\right.
\end{eqnarray}

\begin{definition}
A pair
$(Z,
\eta)$ is called a solution to the problem (\ref{2.9}) if it satisfies\\
(1). $Z=(Z_1,...,Z_{n-1})\in C([0,\infty )\rightarrow R^{n-1})$ and $Z(t)\geq -V(t),$\\
(2). $\eta=(\eta_1,...,\eta_{n-1})\in C([0,\infty )\rightarrow R^{n-1})$ and for each $i$, $\eta_i(t)$ is an increasing continuous
function with
$\eta_i(0)=0$,\\
(3). for all $t\geq 0$,
$$Z(t)=\int_0^tn^2A^nZ(s)ds+\eta(t),$$\\
(4). for $t\geq 0$, $\sum_{i=1}^{n-1}\int_0^t
(Z_i(s)+V_i(s))\eta_i(ds)=0.$
\end{definition}
To prove the existence of the solution to equation (\ref{2.9}), we need the following  estimate which also plays an important role in the subsequent sections.
\begin{lemma}
If $(Z^i(t), \eta^i(t))$ is a solution to equation (\ref{2.9}) with $V$ replaced by $V^i$, $i=1,2$, then for $k\geq 1$, $T>0$, we have
\begin{equation}\label{2.10}
\sup_{0\leq t\leq T}|Z_k^1(t)-Z_k^2(t)|\leq \sup_{0\leq t\leq T, 1\leq j\leq n-1}|V^1_j(t)-V_j^2(t)|.
\end{equation}
\end{lemma}
\vskip 0.4cm
\noindent{\bf Proof}. Set $m:=\sup_{0\leq t\leq T, 1\leq j\leq n-1}|V^1_j(t)-V_j^2(t)|$ and $M=(m,m,...,m)\in R^{n-1}$.
From the definition of the matrix $A^n$, it is easy to see  that $A^nM=(-m,0,...0,-m)$. Thus we have
$$d(Z^1(t)-Z^2(t)-M)=n^2A^n(Z^1(t)-Z^2(t)-M)dt +n^2A^nMdt+d\eta^1(t)-d\eta^2(t).$$
By the chain rule,
\begin{eqnarray}\label{2.11}
&&d\left [\sum_{k=1}^{n-1}\left ( (Z^1_k(t)-Z^2_k(t)-m )^+\right )^2\right ]\nonumber\\
&=&2\sum_{k=1}^{n-1}(Z^1_k(t)-Z^2_k(t)-m )^+d(Z^1_k(t)-Z^2_k(t)-m )\nonumber\\
&=&2n^2\left [\sum_{k=1}^{n-1}(Z^1_k(t)-Z^2_k(t)-m )^+\sum_{i=1}^{n-1}A^n_{ki}(Z^1_i(t)-Z^2_i(t)-m )\right ]dt\nonumber\\
&&+2n^2<(Z^1(t)-Z^2(t)-M)^+, A^nM>dt\nonumber\\
&&+2\sum_{k=1}^{n-1}(Z^1_k(t)-Z^2_k(t)-m )^+d\eta^1_k(t)\nonumber\\
&&-2\sum_{k=1}^{n-1}(Z^1_k(t)-Z^2_k(t)-m )^+d\eta^2_k(t)\nonumber\\
&:=&I_1+I_2+I_3+I_4.
\end{eqnarray}
By Lemma 3.1,
\begin{equation}\label{2.12}
I_1=2n^2<(Z^1(t)-Z^2(t)-M)^+, A^n(Z^1(t)-Z^2(t)-M)>\leq 0.
\end{equation}
In view of the expression of $A^nM$, we have
\begin{equation}\label{2.12-1}
I_2=2n^2[-m(Z^1_1(t)-Z^2_1(t)-m )^+-m(Z^1_{n-1}(t)-Z^2_{n-1}(t)-m )^+]\leq 0.
\end{equation}
Observe that
\begin{eqnarray*}
&&\{t; Z^1_k(t)-Z^2_k(t)>m\}\nonumber\\
&\subset& \{t; Z^1_k(t)>Z^2_k(t)+m\}\nonumber\\
&\subset& \{t; Z^1_k(t)>-V^2_k(t)+m\}\nonumber\\
&\subset& \{t; Z^1_k(t)>-V^1_k(t)\}.
\end{eqnarray*}
Therefore,
\begin{eqnarray}\label{2.13}
I_3&\leq & 2\sum_{k=1}^{n-1}(Z^1_k(t)-Z^2_k(t)-m )^+\chi_{\{t; Z^1_k(t)>-V_k^1(t)\}}d\eta^1_k(t)\nonumber\\
&=& 0.
\end{eqnarray}
Clearly, $I_4\leq 0$ because of the negative sign.
It follows from (\ref{2.11})--(\ref{2.13}) that
$$d\left [\sum_{k=1}^{n-1}\left ( (Z^1_k(t)-Z^2_k(t)-m )^+\right )^2\right ]\leq 0.$$
Hence
$$\sum_{k=1}^{n-1}\left ( (Z^1_k(t)-Z^2_k(t)-m)^+\right )^2\leq \sum_{k=1}^{n-1}\left ((-m )^+\right )^2=0$$
proving the Lemma. $\square$

\begin{theorem}
There exists a unique solution $(Z, \eta)$  to the system (\ref{2.9}).
\end{theorem}
\noindent{\bf Proof}.
We first prove the existence. Assume for the moment $V\in C^1([0, \infty )\rightarrow R^{n-1})$. Consider the following system with reflecting boundary on the convex domain $D_n=\{(z_1,...,z_{n-1}); z_k\geq 0, k=1,...,n-1\}$:
\begin{eqnarray}\label{2.14}\left\{
\begin{array}{l}
du(t)=n^2A^nu(t)dt-n^2A^nV(t)dt+V^{\prime}(t)dt+d\eta(t);\\
u_i(t)\geq 0, \quad i=1,..., n-1;\\
\int_0^tu_i(s)d\eta_i(t)=0, \quad i=1,..., n-1.
\end{array}\right.
\end{eqnarray}
It is well known that the above system admits a unique solution $(u, \eta)$, see \cite{LS}. Let $Z(t):=u(t)-V(t)$. It is easy to verify that $(Z, \eta )$ is the unique solution to the system (\ref{2.9}).
Now consider the general case $V\in C([0, \infty )\rightarrow R^{n-1})$. Take a sequence $V^m\in C^1([0, \infty )\rightarrow R^{n-1})$, $m\geq 1$,
that converges to $V$ uniformly on any finite interval. Let $(Z^m, \eta^m)$ denote the unique solution to the system:
\begin{eqnarray}\label{2.15}\left\{
\begin{array}{l}
dZ^m(t)=n^2A^nZ^m(t)dt+d\eta^m(t);\\
Z^m(t)\geq -V^m(t);\\
\int_0^T<Z^m(t)+V^m(t),d\eta^m(t)>=0.
\end{array}\right.
\end{eqnarray}
By Lemma 3.2 it follows that for $T>0$,
\begin{eqnarray*}
&&\lim_{m,l\rightarrow \infty}\sup_{0\leq t\leq T}|Z^m(t)-Z^l(t)|\nonumber\\
&\leq & \lim_{m,l\rightarrow \infty}sup_{0\leq t\leq T}|V^m(t)-V^l(t)|=0.
\end{eqnarray*}
Thus there exists $Z\in C([0, \infty )\rightarrow R^{n-1})$ such that $Z^m\rightarrow Z$ uniformly on finite intervals. From the equation (\ref{2.15}) we see that $\eta^m$ also converges uniformly on finite intervals to some $\eta\in C([0, \infty )\rightarrow R^{n-1})$. Furthermore, letting $m\rightarrow \infty$ in (\ref{2.15}), we see that $(Z, \eta )$ is a
solution to the system (\ref{2.9}).\\
We show now the uniqueness. Let $(Z, \eta )$, $(\hat{Z}, \hat{\eta})$  be two solutions to the system (\ref{2.9}). By the chain rule,
\begin{eqnarray}\label{2.16}
&&|Z(t)-\hat{Z}(t)|^2\nonumber\\
&=&2n^2\int_0^t<Z(s)-\hat{Z}(s), A^n(Z(s)-\hat{Z}(s))>ds\nonumber\\
&&+2\int_0^t<Z(s)-\hat{Z}(s), d\eta(s)-d\hat{\eta}(s)>\nonumber\\
&\leq &2\int_0^t<Z(s)+V(s)-V(s)-\hat{Z}(s), d\eta(s)-d\hat{\eta}(s)>\nonumber\\
&=&-2\int_0^t<V(s)+\hat{Z}(s), d\eta(s)>-2\int_0^t<Z(s)+V(s),d\hat{\eta}(s)>\nonumber\\
&\leq& 0,
\end{eqnarray}
where we have used the fact that $V(s)+\hat{Z}(s)\geq 0$, $V(s)+Z(s)\geq 0$ (as vectors).
Hence, $Z=\hat{Z}$ which further implies $\eta=\hat{\eta}$ from the equation (\ref{2.9}). $\square$
\section{A discretization scheme for deterministic obstacle problems}
\setcounter{equation}{0}
In this section, we will introduce a discretization scheme for parabolic obstacle problems and establish the convergence of the scheme.
Consider  the following  parabolic obstacle problem:
\begin{eqnarray}\label{3.1-1}\left\{
\begin{array}{l}
\frac{\partial Z(t,x)}{\partial t}-\frac{\partial^2 Z(t,x)}{\partial x^2}=\dot{\eta}(t,x), \quad x\in [0,1];\\
Z(t,x)\geq -V(t,x);\\
\int_0^t\int_0^1(Z(s,x)+V(s,x))\eta(ds,dx)=0,
\end{array}\right.
\end{eqnarray}
where $V\in C(R_+\times [0,1])$ with $V(0,x)=u_0(x)\geq 0$.
\begin{definition}
 If a pair
$(Z,
\eta)$ satisfies\\
(1). $Z$ is a continuous function on $R_+\times [0,1]$ and
$$Z(0,x)=0, Z(t,0)=Z(t,1)=0,\  Z(t,x)\geq -V(t,x),$$
(2). \ $\eta$ is a measure on $(0,1)\times \mathbb{R}_+$ such that
for all $\varepsilon>0, T>0$
$$\eta\big ([0,T]\times (\varepsilon,(1-\varepsilon))\big )<\infty,$$
(3). for all $t\geq 0, \phi\in C_0^2(0,1)$,
\begin{eqnarray*}\label{weak form}
&&(Z(t),\phi)-\int_0^t(Z(s),\phi'')ds
=\int_0^t\int_0^1\phi(x)\eta(ds, dx),\
\end{eqnarray*}
(4). $\int_0^t\int_0^1
(Z(s,x)+V(s,x))\eta(ds,dx)=0, \quad t\geq 0$ \\
then $(Z, \eta)$ is called a solution to problem
(\ref{3.1-1}).
\end{definition}

 The following result was proved
in \cite{NP}(Theorem 1.4).

\begin{proposition} $({\bf [NP]})$
If $V(0,x)=u_0(x), V(t,0)=V(t,1)=0$ for all $ t\geq 0,$ Eq.
(\ref{3.1-1}) admits a unique solution. Moreover, if $Z^1$, $Z^2$ are solutions of the obstacle problem (\ref{3.1-1}) with $V$ replaced respectively
with $V^1$ and $V^2$, then
$|Z^1-Z^2|^T_\infty\leq |V^1-V^2|^T_\infty,$ for $T>0$. Where $|Z^1-Z^2|^T_\infty=\sup_{0\leq t\leq T, 0\leq x\leq 1}|Z^1(t,x)-Z^2(t,x)|$ and $|V^1-V^2|^T_\infty$ is defined accordingly.
\end{proposition}
We now introduce the discretization scheme for the deterministic obstacle problem (\ref{3.1-1}). For very positive integer $n\geq 1$, define
$$V^n(t)=\left(V(t, \frac{1}{n}), ...,V(t, \frac{n-1}{n})\right ),$$
where $V(t,x)$ is the function appeared in equation (\ref{3.1-1}).

Consider  the following Skorohod-type reflecting system in $R^{n-1}$:
\begin{eqnarray}\label{3.1-2}\left\{
\begin{array}{l}
dZ^n(t)=n^2A^nZ^n(t)dt+d\eta^n(t);\\
Z^n(t)\geq -V^n(t);\\
\int_0^T<Z^n(t)+V^n(t),d\eta^n(t)>=0.
\end{array}\right.
\end{eqnarray}
The existence and uniqueness of the solution of the above system was proved in Section 3.
For $n\geq 1$, define the  continuous functions $Z^n$ by
\begin{equation}\label{3.2-3}
Z^n(t,x):=Z_k^n(t)+(nx-k)\left ( Z_{k+1}^n(t)-Z_{k}^n(t)\right )
\end{equation}
for $x\in [\frac{k}{n}, \frac{k+1}{n})$, $k=0,...,n-1$, where $Z_0^n(t), Z_n^n(t)$ are set to be zero.
We have
\begin{theorem}
Let $Z$ be the solution to equation (\ref{3.1-1}). Then for $T>0$,
\begin{equation}\label{3.2-4}
\lim_{n\rightarrow\infty}\sup_{0\leq t\leq T, 0\leq x\leq 1}|Z^n(t,x)-Z(t,x)|=0.
\end{equation}
\end{theorem}
\noindent{\bf Proof}. We divide the proof into two steps.

{\bf Step 1}. Suppose $V\in C^{1,2}([0, \infty )\times [0,1])$. In this case we first show that the function $\eta^n(t)$ in (\ref{3.1-2}) is absolutely continuous and
\begin{equation}\label{3.2-4}
\int_0^T|\dot{\eta}^n|^2(t)dt=\int_0^T\sum_{k=1}^{n-1}(\dot{\eta}_k^n(t))^2dt\leq C(\int_0^T|\dot{V}^n(t)|^2dt + n),
\end{equation}
for some constant $C$ independent of $n$, where $\dot{V}^n(t)$ stands for the derivative of $V^n$.
Indeed, let $U^n(t):=Z^n(t)+V^n(t)$. Then $(U^n, \eta^n)$ is the solution of the reflecting system:
\begin{eqnarray}\label{3.1-5}\left\{
\begin{array}{l}
dU^n(t)=n^2A^nU^n(t)dt+\dot{V}^n(t)dt-n^2A^nV^n(t)dt +d\eta^n(t);\\
U^n(t)\geq 0;\\
\int_0^T<U^n(t),d\eta^n(t)>=0.
\end{array}\right.
\end{eqnarray}
Define $\phi(z):=\sum_{k=1}^{n-1}(z_k^-)^2$ for $z\in R^{n-1}$, where $z_k^-$ stands for the negative part of $z_k$.  Consider the following penalized equation:
\begin{eqnarray}\label{3.1-6}
dU^{n,\varepsilon}(t)&=&n^2A^nU^{n,\varepsilon}(t)dt+\dot{V}^n(t)dt-n^2A^nV^n(t)dt\nonumber\\
&&-\frac{1}{\varepsilon} \nabla \phi (U^{n,\varepsilon}(t))dt.
\end{eqnarray}
According to \cite{LS}, it holds that
\begin{equation}\label{3.1-7}
\lim_{\varepsilon\rightarrow 0}\sup_{0\leq t\leq T}|U^{n,\varepsilon}(t)-U^n(t)|=0,
\end{equation}
\begin{equation}\label{3.1-8}
\eta^n(t)=-\lim_{\varepsilon\rightarrow 0}\frac{1}{\varepsilon} \int_0^t\nabla \phi (U^{n,\varepsilon}(s))ds, \quad\quad \mbox{for}\quad t>0.
\end{equation}
Using the chain rule we have
\begin{eqnarray}\label{3.1-9}
\phi(U^{n,\varepsilon}(t))&=&n^2\int_0^t<\nabla\phi (U^{n,\varepsilon}(s)), A^nU^{n,\varepsilon}(s)>ds+\int_0^t<\nabla\phi (U^{n,\varepsilon}(s)), \dot{V}^n(s)>ds\nonumber\\
&&-n^2\int_0^t<\nabla\phi (U^{n,\varepsilon}(s)), A^nV^n(s)>ds-\frac{1}{\varepsilon} \int_0^t|\nabla \phi|^2 (U^{n,\varepsilon}(s))ds.
\end{eqnarray}
As in the proof of  lemma 3.2 , we can show that $
<b^-, A^nb>\geq 0$  for all $ b\in R^{n-1}$. Thus
\begin{eqnarray}\label{3.1-10}
&&n^2\int_0^t<\nabla\phi (U^{n,\varepsilon}(s)), A^nU^{n,\varepsilon}(s)>ds\nonumber\\
&=&-2n^2\int_0^t<(U^{n,\varepsilon}(s))^-, A^nU^{n,\varepsilon}(s)>ds\leq 0.
\end{eqnarray}
As $\phi \geq 0$, it follows from (\ref{3.1-9}) and (\ref{3.1-10}) that
\begin{eqnarray}\label{3.1-11}
&&\frac{1}{\varepsilon} \int_0^t|\nabla \phi|^2 (U^{n,\varepsilon}(s))ds\nonumber\\
&\leq &\int_0^t<\nabla\phi (U^{n,\varepsilon}(s)), \dot{V}^n(s)>ds-n^2\int_0^t<\nabla\phi (U^{n,\varepsilon}(s)), A^nV^n(s)>ds\nonumber\\
&\leq& \left (\int_0^t|\nabla \phi|^2 (U^{n,\varepsilon}(s))ds\right )^{\frac{1}{2}}\times\{(\int_0^t|\dot{V}^n(s)|^2ds)^{\frac{1}{2}}
+(\int_0^t|n^2A^nV^n(s)|^2ds)^{\frac{1}{2}}\}.
\end{eqnarray}
which yields that
\begin{eqnarray}\label{3.1-12}
&& \int_0^t|\frac{1}{\varepsilon} \nabla \phi|^2 (U^{n,\varepsilon}(s))ds\nonumber\\
&\leq &C\{\int_0^t|\dot{V}^n(s)|^2ds
+\int_0^t|n^2A^nV^n(s)|^2ds\}, \quad\quad \mbox{for all} \quad t>0.
\end{eqnarray}
By selecting a subsequence if necessary, we conclude that $\frac{1}{\varepsilon} \nabla \phi (U^{n,\varepsilon}(\cdot))$ converges weakly in $L^2([0, T] \rightarrow R^{n-1})$ as $\varepsilon \rightarrow 0$. Combing with (\ref{3.1-8}) we deduce that $\eta^n(t)$ is absolutely continuous and
\begin{eqnarray}\label{3.1-13}
\int_0^T|\dot{\eta}^n|^2(t)dt&\leq& \liminf_{\varepsilon\rightarrow 0}\int_0^T|\frac{1}{\varepsilon} \nabla \phi|^2 (U^{n,\varepsilon}(s))ds\nonumber\\
&\leq &C\{\int_0^T|\dot{V}^n(s)|^2ds
+\int_0^T|n^2A^nV^n(s)|^2ds\}, \quad\quad \mbox{for all} \quad t>0.
\end{eqnarray}
From the definition of $A^n$, it is seen that
\begin{equation}\label{3.1-14}
A^nV^n(t)=\left(\begin{array}{l}a_1\\  a_2\\ \cdot \\\cdot\\\cdot\\a_{n-1}\end{array}\right ), where \quad a_k=V(t, \frac{k+1}{n})-2V(t, \frac{k}{n})+V(t, \frac{k-1}{n}).
\end{equation}
Observe that
\begin{eqnarray}\label{3.1-15}
|n^2a_k|&=&n^2\left |\int_{\frac{k}{n}}^{\frac{k+1}{n}}dy\int_{\frac{k}{n}}^y\frac{\partial^2V(t,z)}{\partial z^2}dz+\int_{\frac{k-1}{n}}^{\frac{k}{n}}dy\int_y^{\frac{k}{n}}\frac{\partial^2V(t,z)}{\partial z^2}dz \right |\nonumber\\
&\leq &2\sup_{0\leq t\leq T, 0\leq z\leq 1}|\frac{\partial^2V(t,z)}{\partial z^2}|.
\end{eqnarray}
Substitute (\ref{3.1-15}) back to (\ref{3.1-13}) to complete the proof of (\ref{3.2-4}). Next we show that the family  $\{Z^n(t,x), n\geq 1\}$ defined in (\ref{3.2-3}) is relatively compact in the space
$C([0,T]\times [0,1])$.  Recall $G^n(t)=e^{n^2A^nt}$ as in Section 2. By the variation of constant formula, we have
\begin{equation}\label{3.1-16}
Z^n(t)=G^n(t)u^n(0)+\int_0^tG^n(t-s)\dot{\eta}^n(s)ds.
\end{equation}
For $n\geq 1$, define
\begin{equation}\label{3.1-17}
\dot{\eta}^n(t,x)=\dot{\eta}^n_k(t)+(nx-k)(\dot{\eta}^n_{k+1}(t)-\dot{\eta}^n_k(t))\quad \mbox{for} \quad x\in [\frac{k}{n}, \frac{k+1}{n}), \quad k=0,...,n-1,
\end{equation}
with $\dot{\eta}^n_0(t):=0, \dot{\eta}^n_n(t):=0$.
Set $\varphi_j(x):=\sqrt{2}sin(jx\pi)$. As in \cite{G} introduce the kernel $G^n(t,x,y)$ by
\begin{equation}\label{3.1-17-0}
G^n(t,x,y)=\sum_{j=1}^{n-1}exp(\lambda_j^nt)\varphi_j^n(x)\varphi_j(k_n(y)),
\end{equation}
where $k_n(y)=\frac{[ny]}{n}$ and for $x\in [\frac{k}{n}, \frac{k+1}{n}]$, define
\begin{equation}\label{3.1-17-1}
\varphi_j^n(x)=\varphi_j(\frac{k}{n})+(nx-k)(\varphi_j(\frac{k+1}{n})-\varphi_j(\frac{k}{n})).
\end{equation}
The following statements were proved in \cite{G}(see the proof of Lemma 3.6 there).
\begin{equation}\label{3.1-18}
\int_0^{s}\int_0^1|G^n(t-r,x,y)-G^n(s-r,x,y)|^2drdy\leq C_1\sqrt{t-s},
\end{equation}
for $x\in [0,1]$ and $s\leq t\leq T$.

\begin{equation}\label{3.1-19}
\int_s^{t}\int_0^1|G^n(t-r,x,y)|^2drdy\leq C_2\sqrt{t-s},
\end{equation}
for $x\in [0,1]$ and $s\leq t\leq T$.

\begin{equation}\label{3.1-20}
\int_0^{t}\int_0^1|G^n(t-r,x,z)-G^n(t-r,y,z)|^2drdz\leq C_3|x-y|,
\end{equation}
for $x, y\in [0,1]$ and $0\leq t\leq T$.\\
The constants $C_1$, $C_2$, $C_3$ in the above estimates  are independent of $n$.
\vskip 0.3cm
By (\ref{3.1-16}) and  a simple calculation we find that
\begin{equation}\label{3.1-21}
Z^n(t,x)=\int_0^t\int_0^1G^n(t-s,x,y)\dot{\eta}^n(s,k_n(y))dsdy.
\end{equation}

The estimate (\ref{3.2-4}) yields that
\begin{eqnarray}\label{3.1-22}
&&\int_0^T\int_0^1|\dot{\eta}^n|^2(s,k_n(y)))dsdy
=\sum_{k=0}^{n-1}\frac{1}{n}\int_0^T|\dot{\eta}_k^n|^2(s)ds\nonumber\\
&\leq&C(\int_0^T\int_0^1\frac{\partial V(t, k_n(y))}{\partial t}|^2dtdy + 1)\nonumber\\
&\leq& C_T,
\end{eqnarray}
where we have used the smoothness assumptions on $V$ and the definition of $V^n$. Using the above estimate and H\"o{}lder's inequality  it follows easily from (\ref{3.1-21}), (\ref{3.1-20}), (\ref{3.1-19})
and (\ref{3.1-18}) that there exists a constant $C$, independent of $n$, such that
\begin{equation}\label{3.1-23}
|Z^n(t,x)-Z^n(s,y)|^2\leq C(\sqrt{|t-s|}+|x-y|), \quad\quad s,t\in [0,T], x,y\in [0,1].
\end{equation}
By Arzela-Ascoli theorem, $\{ Z^n(t,x), n\geq 1\}$ is relatively compact. On the other hand, (\ref{3.1-22}) implies that $\{ \dot{\eta}^n(\cdot,k_n(\cdot)), n\geq 1\}$ is relatively compact in $L^2([0,T]\times [0,1])$ with respect to the weak topology.
Selecting a subsequence if necessary, we can assume that $Z^n(\cdot, \cdot)$ converges uniformly to some function $Z(\cdot,\cdot) \in C([0,T]\times [0,1])$ and $\dot{\eta}^n(\cdot,k_n(\cdot))$ converges weakly to some $\dot{\eta}(\cdot,\cdot)\in L^2([0,T]\times [0,1])$. We complete the proof of step 1 by showing that $(Z, \eta(dt,dy):=\dot{\eta}(t,y)dtdy)$ is the solution to the system (\ref{3.1-1}). For  $\phi\in
C^2_0((0,1))$, set $\phi^n:=(\phi(\frac{1}{n}),...,\phi(\frac{n-1}{n}))$. By the symmetry of the matrix $A^n$ it follows from (\ref{3.2-4}) that
\begin{eqnarray}\label{3.1-24}
<Z^n(t), \phi^n>&=&\int_0^t<n^2A^n\phi^n, Z^n(s)>ds+\int_0^t<\phi^n, \dot{\eta}^n(s)>ds.
\end{eqnarray}
Multiply the above equation by $\frac{1}{n}$ to get
\begin{eqnarray}\label{3.1-25}
\int_0^1 Z^n(t, k_n(y))\phi(k_n(y))dy&=&\int_0^tds\int_0^1\Delta_n\phi(k_n(y)) Z^n(s,k_n(y))dy\nonumber\\
&&+\int_0^tds\int_0^1\phi(k_n(y)) \dot{\eta}^n(s,k_n(y))dy,
\end{eqnarray}
where $\Delta_n\phi(x):=n^2(\phi(x+\frac{1}{n})-2\phi(x)+\phi(x-\frac{1}{n}))$ is the discrete Laplacian operator.
Letting $n\rightarrow \infty$ in (\ref{3.1-25})  we obtain
\begin{eqnarray}\label{3.1-26}
\int_0^1 Z(t,y)\phi(y)dy&=&\int_0^tds\int_0^1\phi^{\prime\prime}(y) Z(s,y)dy+\int_0^tds\int_0^1\phi(y) \dot{\eta}(s,y)dy,
\end{eqnarray}
where we have used the fact that $\phi(k_n(y)) \rightarrow \phi(y)$ (strongly) in $L^2([0,T]\times [0,1])$.
On the other hand, it follows from the definition that
\begin{equation}\label{3.1-26-0}
\int_0^t\int_0^1(Z^n(s, k_n(y))+V(s, k_n(y)))\dot{\eta}^n(s,k_n(y))dsdy=0.
\end{equation}
Invoking (\ref{3.1-23}) and the dominated convergence theorem we have
\begin{eqnarray}\label{3.1-26-1}
&&\int_0^T\int_0^1 (Z^n(s, k_n(y))+V(s, k_n(y))-  Z(s,y)-V(s,y))^2dsdy\nonumber\\
&\leq& C\int_0^T\int_0^1 (Z^n(s, k_n(y))-Z^n(s,y))^2dsdy+C\int_0^T\int_0^1 (Z^n(s,y)-Z(s,y))^2dsdy\nonumber\\
&&\quad \quad +\int_0^T\int_0^1 (V(s, k_n(y))-V(s,y))^2dsdy\nonumber\\
&\leq& C(\frac{1}{n})^2 +C\int_0^T\int_0^1 (Z^n(s,y)-Z(s,y))^2dsdy+\int_0^T\int_0^1 (V(s, k_n(y))-V(s,y))^2dsdy\nonumber\\
&&\rightarrow 0 \quad\quad\quad \mbox{as}\quad \quad n\rightarrow \infty.
\end{eqnarray}
Letting $n\rightarrow \infty$ in {\ref{3.1-26-0}), the weak convergence of $\dot{\eta}^n$ and (\ref{3.1-26-1}) yield 
$$\int_0^t\int_0^1(Z(s, y)+V(s, y))\dot{\eta}(s,y)dsdy=0.$$
We have shown that  conditions (3), (4) in the Definition 4.1 are satisfied by $(Z,\eta)$. It is straightforward to also check that $(Z,\eta)$ satisfies (1),(2) in the definition 4.1.  Thus, $(Z,\eta)$ is the solution to equation (\ref{3.1-1}).
\vskip 0.4cm
{\bf Step 2}. The general case $V\in C(R_+\times [0,1])$.
Take a sequence $V_m\in C^{1,2}(R_+\times [0,1]), m\geq 1$ such that $\sup_{0\leq t\leq T, 0\leq x\leq 1}|V_m(t,x)-V(t,x)|\rightarrow 0$ as $m\rightarrow \infty$ for any $T>0$. For every integer $n\geq 1$, define
$$V^{m,n}(t)=\left(V_m(t, \frac{1}{n}), ...,V_m(t, \frac{n-1}{n})\right ).$$
Let $(Z^{m,n}, \eta^{m,n})$ be the solution to the following Skorohod-type problem in $R^{n-1}$:
\begin{eqnarray}\label{3.1-27}\left\{
\begin{array}{l}
dZ^{m,n}(t)=n^2A^nZ^{m,n}(t)dt+d\eta^{m,n}(t);\\
Z^{m,n}(t)\geq -V^{m,n}(t);\\
\int_0^T<Z^{m,n}(t)+V^{m,n}(t),d\eta^{m,n}(t)>=0.
\end{array}\right.
\end{eqnarray}
Set $Z_0^{m,n}(t)=0$ and introduce the  continuous functions $Z^{m,n}$ by
\begin{equation}\label{3.1-28}
Z^{m,n}(t,x):=Z_k^{m,n}(t)+(nx-k)\left ( Z_{k+1}^{m,n}(t)-Z_{k}^{m,n}(t)\right )
\end{equation}
for $x\in [\frac{k}{n}, \frac{k+1}{n})$, $k=0,...,n-1$.
According to the result proved in step 1, for every $m\geq 1$ we have
\begin{equation}\label{3.1-29}
\lim_{n\rightarrow\infty}\sup_{0\leq t\leq T, 0\leq x\leq 1}|Z^{m,n}(t,x)-Z^{(m)}(t,x)|=0,
\end{equation}
where $Z^{(m)}(t,x)$ is the solution of   the following  parabolic obstacle problem:
\begin{eqnarray}\label{3.1-30}\left\{
\begin{array}{l}
\frac{\partial Z^{(m)}(t,x)}{\partial t}-\frac{\partial^2 Z^{(m)}(t,x)}{\partial x^2}=\dot{\eta}^{(m)}(t,x);\\
Z^{(m)}(t,x)\geq -V_m(t,x);\\
\int_0^t\int_0^1(Z^{(m)}(t,x)+V_m(s,x))\eta^{(m)}(ds,dx)=0.
\end{array}\right.
\end{eqnarray}
Applying Lemma 3.2 we have
\begin{eqnarray}\label{3.1-31}
&&\sup_{0\leq t\leq T, 0\leq x\leq 1}|Z^{m,n}(t,x)-Z^n(t,x)|\nonumber\\
&=&\sup_{0\leq t\leq T, 1\leq k\leq n-1}|Z^{m,n}_k(t)-Z^n_k(t)|\nonumber\\
&\leq &\sup_{0\leq t\leq T, 1\leq k\leq n-1}|V_k^{m,n}(t)-V_k^n(t)|\nonumber\\
&=&\sup_{0\leq t\leq T, 0\leq x\leq 1}|V_m(t,k_n(x))-V(t,k_n(x))|\nonumber\\
&\leq& \sup_{0\leq t\leq T, 0\leq x\leq 1}|V_m(t,x)-V(t,x)|.
\end{eqnarray}
Now we are in the position to complete the proof of the theorem. For every $m\geq 1$, by (\ref{3.1-31}) and Proposition 4.1  we have
\begin{eqnarray}\label{3.1-32}
&&\sup_{0\leq t\leq T, 0\leq x\leq 1}|Z^{n}(t,x)-Z(t,x)|\nonumber\\
&\leq&\sup_{0\leq t\leq T, 0\leq x\leq 1}|Z^{n}(t,x)-Z^{(m)}(t,x)|+\sup_{0\leq t\leq T, 0\leq x\leq 1}|Z^{(m)}(t,x)-Z(t,x)|\nonumber\\
&\leq&\sup_{0\leq t\leq T, 0\leq x\leq 1}|Z^{n}(t,x)-Z^{m,n}(t,x)|+\sup_{0\leq t\leq T, 0\leq x\leq 1}|Z^{m,n}(t,x)-Z^{(m)}(t,x)|\nonumber\\
&&+\sup_{0\leq t\leq T, 0\leq x\leq 1}|V_m(t,x)-V(t,x)|\nonumber\\
&\leq& 2\sup_{0\leq t\leq T, 0\leq x\leq 1}|V_m(t,x)-V(t,x)|+\sup_{0\leq t\leq T, 0\leq x\leq 1}|Z^{m,n}(t,x)-Z^{(m)}(t,x)|.
\end{eqnarray}
Given a positive constant $\varepsilon>0$. First choose $m$ sufficiently large such that
\begin{equation}\label{3.1-33}
 2\sup_{0\leq t\leq T, 0\leq x\leq 1}|V_m(t,x)-V(t,x)|\leq \frac{\varepsilon}{2}.
 \end{equation}
For such a fixed $m$, applying (\ref{3.1-29}) there exists an integer $N$ such that for $n\geq N$,
\begin{equation}\label{3.1-34}
 \sup_{0\leq t\leq T, 0\leq x\leq 1}|Z^{m,n}(t,x)-Z^{(m)}(t,x)|\leq \frac{\varepsilon}{2}.
 \end{equation}
 Putting (\ref{3.1-32}), (\ref{3.1-33}) and (\ref{3.1-34}) together we obtain that
 $$\sup_{0\leq t\leq T, 0\leq x\leq 1}|Z^{n}(t,x)-Z(t,x)|\leq \varepsilon$$
 for $n\geq N$.
 As $\varepsilon$ is arbitrary, the proof is complete.$\square$

\section{The convergence of the scheme}
\setcounter{equation}{0}
After all the preparations in the previous sections, this part is devoted to the proof of the  main result.
For $y\in R^{n-1}$, set
$$F_n(t,y)=\left (\begin{array}{l}
f(t,\frac{1}{n},y_1)\\
f(t,\frac{2}{n},y_2)\\
.\\
.\\
.\\
f(t,\frac{n-1}{n},y_{n-1})
\end{array}
\right )
$$
$$
\Sigma_n(t,y)=\left (\begin{array}{llll}
\sigma(t,\frac{1}{n},y_1)&0& ..&0\\
0& \sigma(t,\frac{2}{n},y_2)&..&0\\
.&.&.&.\\
.&.&.&.\\
.&.&.&.\\
0&0&..&\sigma(t,\frac{n-1}{n},y_{n-1})
\end{array}
\right )
$$
The system (\ref{2.4}) can be written as
\begin{eqnarray}\label{3.2}
du^n(t)&=&n^2A^nu^n(t)dt+F_n(t, u^n(t))dt\nonumber\\
&+& \sqrt{n}\Sigma_n(t,u^n(t))dW^n(t)+d\eta^n(t),
\end{eqnarray}
where $W^n=(W^n_1(t),...,W^n_{n-1}(t))$.

By the variation of constant formula, it follows that
\begin{eqnarray}\label{3.3}
u^n(t)&=&G^n(t)u^n(0)+\int_0^tG^n(t-s)F_n(s, u^n(s))ds\nonumber\\
&+& \sqrt{n}\int_0^tG^n(t-s)\Sigma_n(s,u^n(s))dW^n(s)+\int_0^tG^n(t-s)d\eta^n(s),
\end{eqnarray}
where as before $G^n(t)=exp(n^2A^nt)$.
Denote
\begin{eqnarray}\label{3.3-1}
v^n(t)&=&G^n(t)u^n(0)+\int_0^tG^n(t-s)F_n(s, u^n(s))ds\nonumber\\
&&+\sqrt{n}\int_0^tG^n(t-s)\Sigma_n(s,u^n(s))dW^n(s).
\end{eqnarray}
Then $v^n$ is the solution of the SDE:
\begin{eqnarray}\label{3.4}
dv^n(t)&=&n^2A^nv^n(t)dt+F_n(t, u^n(t))dt\nonumber\\
&+& \sqrt{n}\Sigma_n(t,u^n(t))dW^n(t),
\end{eqnarray}
and $(Z^n(t):=u^n(t)-v^n(t), \eta^n(t))$ is the solution of the system:
\begin{eqnarray}\label{3.4-1}\left\{
\begin{array}{l}
dZ^n(t)=n^2A^nZ^n(t)dt+d\eta^n(t);\\
Z^n(t)\geq -v^n(t);\\
\int_0^T<Z^n(t)+v^n(t),d\eta^n(t)>=0, \quad\quad \mbox{for all} \quad T>0.
\end{array}\right.
\end{eqnarray}
Recall the random field $u^n(t,x)$ defined in (\ref{2.6}) in Section 2 and introduce the random fields
\begin{equation}\label{3.9}
\eta^n(t,x)=\eta^n_k(t)+(nx-k)(\eta^n_{k+1}(t)-\eta^n_k(t))\quad \mbox{for} \quad x\in [\frac{k}{n}, \frac{k+1}{n}), \quad k=0,...,n-1,
\end{equation}
\begin{equation}\label{3.10}
v^n(t,x)=v^n_k(t)+(nx-k)(v^n_{k+1}(t)-v^n_k(t))\quad \mbox{for} \quad x\in [\frac{k}{n}, \frac{k+1}{n}), \quad k=0,...,n-1,
\end{equation}
where $\eta_0^n(t):=0, \eta_n^n(t):=0$, $v_0^n(t):=0, v_n^n(t):=0$.
Let the kernel $G^n(t,x, y)$ be defined as in (\ref{3.1-17-0}) in Section 4. It is easy to verify that $u^n$ and $v^n$ satisfies the  equations
\begin{eqnarray}\label{3.11}
u^n(t,x)&=&\int_0^1G^n(t,x,y)u^n(0,k_n(y))dy\nonumber\\
&&+\int_0^t\int_0^1G^n(t-s,x,y)f(s,k_n(y),u^n(s,k_n(y)))dyds\nonumber\\
&&+\int_0^t\int_0^1G^n(t-s,x,y)\sigma(s,k_n(y),u^n(s,k_n(y)))W(ds,dy)\nonumber\\
&&+\int_0^t\int_0^1G^n(t-s,x,y)\eta^n(ds,k_n(y)))dy,
\end{eqnarray}
\begin{eqnarray}\label{3.12}
v^n(t,x)&=&\int_0^1G^n(t,x,y)u^n(0,k_n(y))dy\nonumber\\
&&+\int_0^t\int_0^1G^n(t-s,x,y)f(s,k_n(y),u^n(s,k_n(y)))dyds\nonumber\\
&&+\int_0^t\int_0^1G^n(t-s,x,y)\sigma(s,k_n(y),u^n(s,k_n(y)))W(ds,dy),
\end{eqnarray}
where $k_n(y)=\frac{[ny]}{n}$ as in Section 4.
Let $G(t,x,y)$ denote the heat kernel of the Laplacian $\frac{\partial^2}{\partial x^2}$ on the interval  $[0,1]$ with the Dirichlet boundary condition, i.e.,
$$G(t,x,y)=\sum_{k=1}^{\infty}exp(-k^2\pi^2t)\varphi_k(x)\varphi_k(y).$$
The following lemma was proved in \cite{G}
\begin{lemma}
The following statements hold:

(i) There exists a constant $c$ such that
\begin{equation}\label{3.13}
\int_0^{\infty}\int_0^1|G(t,x,y)-G^n(t,x,y)|^2dydt\leq \frac{c}{n},
\end{equation}
for $x\in [0,1]$ and $n\geq 1$.

(ii) For every $t>0$, $\gamma\in (0,1)$, $\beta >\frac{\gamma}{2}+\frac{1}{2}$ there is a constant $C$ such that
\begin{equation}\label{3.14}
\int_0^1|G(t,x,y)-G^n(t,x,y)|^2dy\leq  Cn^{-\gamma}t^{-\beta},
\end{equation}
for $x\in [0,1]$ and $n\geq 1$.
\end{lemma}
We have the following estimate for $u^n$.
\begin{lemma}
For any $T>0$, we have
\begin{equation}\label{3.15}
\sup_{0\leq t\leq T, 0\leq x\leq 1}|u^n(t,x)|\leq 2\sup_{0\leq t\leq T, 0\leq x\leq 1}|v^n(t,x)|
\end{equation}
\end{lemma}
\vskip 0.4cm
\noindent{\bf Proof}. Keeping in mind that  $u^n(t,x)$, $v^n(t,x)$ are piecewise linear in $x$, by Lemma 3.2, we have
\begin{eqnarray*}
&&\sup_{0\leq t\leq T, 0\leq x\leq 1}|u^n(t,x)|\nonumber\\
&=&\sup_{0\leq t\leq T, 1\leq k\leq n-1}|u^n(t, \frac{k}{n})|=\sup_{0\leq t\leq T, 1\leq k\leq n-1}|Z^n_k(t)+v^n(t, \frac{k}{n})|\nonumber\\
&\leq& 2\sup_{0\leq t\leq T, 1\leq k\leq n-1}|v^n(t, \frac{k}{n})|=2\sup_{0\leq t\leq T, 0\leq x\leq 1}|v^n(t, x)|
\end{eqnarray*}
proving the lemma. $\square$.

\begin{proposition}
Assume the linear growth condition (H.2) in Section 2. Then for $p\geq 1$ and $T>0$, there exists a constant $C_p$ such that
\begin{equation}\label{3.16}
\sup_nE[\sup_{0\leq t\leq T, 0\leq x\leq 1}|u^n(t,x)|^p]\leq C_p
\end{equation}
\end{proposition}

\noindent {\bf Proof}. We will use the  notation $|u|_{\infty}^t:=\sup_{0\leq s\leq t, 0\leq x\leq 1}|u(s,x)|$. We can assume $p>20$.
By Lemma 5.2, we have
\begin{eqnarray}\label{3.17}
(|u^n|^T_\infty)^p&\leq& 2^p(|v^n|^T_\infty)^p\nonumber\\
&\leq&c(p)(\sup_{x\in [0,1], t\in [0,T]}|\int_0^t\int_0^1G^n(t-s,x,y)f(s,k_n(y), u^n(s,k_n(y)))dyds|)^p\nonumber\\
&&+c(p)(\sup_{x\in [0,1], t\in [0,T]}|\int_0^t\int_0^1G^n(t-s,x,y)\sigma(s,k_n(y), u^n(s,k_n(y)))W(ds,dy)|)^p.\nonumber\\
\end{eqnarray}
Set
$$I_1(t,x):=\int_0^t\int_0^1G^n(t-s,x,y)f(s,k_n(y), u^n(s,k_n(y)))dyds,$$

$$I_2(t,x):=\int_0^t\int_0^1G^n(t-s,x,y)\sigma(s,k_n(y), u^n(s,k_n(y)))W(ds,dy).$$
By the linear growth of $f$ and the  H\"o{}lder inequality,
\begin{eqnarray}\label{3.18}
E(|I_1|^T_\infty)^p&\leq&C_p(T)(\sup_{x\in [0,1], t\in [0,T]}\int_0^t\int_0^1
(G^n(s,x,y))^2dyds )^\frac{p}{2}E\int_0^T(|f(\cdot,k_n(\cdot), u^n(\cdot,k_n(\cdot)))|^t_\infty)^pdt\nonumber\\
&\leq& C_p(T)(1+E\int_0^T(|u^n(\cdot,k_n(\cdot))|^t_\infty)^pdt\nonumber\\
&\leq& C_p(T)(1+E\int_0^T(|u^n|^t_\infty)^pdt
\end{eqnarray}
where $\frac{1}{q}+\frac{1}{p}=1$, $C_p(T)$ denotes a generic constant depending on $T$, $p$. The fact
\begin{equation}\label{3.19}
\sup_n\sup_{x\in [0,1], t\in [0,T]}\int_0^t\int_0^1
(G^n(s,x,y))^2dyds <\infty.
\end{equation}
has also been used in the derivation of (\ref{3.18}). Actually (\ref{3.19}) follows from (\ref{3.13}) and the fact
$$ \sup_{x\in [0,1], t\in [0,T]}\int_0^t\int_0^1
(G(s,x,y))^2dyds <\infty.$$
 In view of (\ref{3.1-18}), (\ref{3.1-19}), (\ref{3.1-20}),
following a similar calculation as in the proof of Corollary
3.4 in \cite{WA} and Lemma 3.6 in \cite{G}, we obtain that
\begin{eqnarray}\label{3.19-1}
&&E|I_2(t,x)-I_2(s,y)|^p\nonumber\\
&\leq&c\big (E\int_0^{(t\vee
s)}(|u^n(\cdot,k_n(\cdot))|^r_\infty)^pdr\big)\times
|(t,x)-(s,y)|^{\frac{p}{4}-3}\nonumber\\
&\leq&c\big (E\int_0^{(t\vee
s)}(|u^n|^r_\infty)^pdr\big)\times
|(t,x)-(s,y)|^{\frac{p}{4}-3}.
\end{eqnarray}
Applying Garsia-Rodemich-Rumsey's lemma (See e.g Theorem 1.1 and Corollary 1.2 in \cite{WA}), we obtain that
\begin{eqnarray}\label{Garsian lemma1}
&&|I_2(t,x)-I_2(s,y)|^p\nonumber\\
&\leq&
N(\omega)^p|(t,x)-(s,y)|^{\frac{p}{4}-5}\big(\log\big(\frac{\gamma}{|(t,x)-(s,y)|}\big)\big)^2,
\end{eqnarray}
where $N(\omega)$ is a random variable satisfying
\begin{eqnarray}\label{estimate on random variable1}
E[N^p]\leq ac\big(E\int_0^{(t\vee
s)}(|u^n|^r_\infty)^pdr\big),
\end{eqnarray}
where $a, \gamma$ are constants depending only on $p$ and $c$ is the constant appeared in (\ref{3.19-1}). Choosing
$s=0$ in $(\ref{Garsian lemma1})$, we see that there exists a
constant $c_T$ such that
\begin{eqnarray}\label{I_2}
E(\sup_{x\in [0,1], t\in [0,T]}|I_2(t,x)|^p) &\leq&c_T E\int_0^T(|u^n|^t_\infty)^pdt.
\end{eqnarray}
Putting $(\ref{3.17}), (\ref{3.18})$, $(\ref{I_2})$ together, we get that
\begin{eqnarray}\label{successive}
E(|u^n|^T_\infty)^p&\leq& c(p,K,T)E\int_0^T(|u^n|^t_\infty)^pdt,
\end{eqnarray}
where $c(p,K,T)$ is a constant depending on $p,K,T.$
Applying the Grownwall's lemma, we proves the proposition. $\square$
 \vskip 0.4cm
 \noindent{\bf Proof of the main result (Theorem 2.1)}.
\vskip 0.3cm
Recall
\begin{eqnarray}\label{3.21}
u(t,x)&=&\int_0^1G(t,x,y)u(0,y)dy\nonumber\\
&&+\int_0^t\int_0^1G(t-s,x,y)f(s,y,u(s,y))dyds\nonumber\\
&&+\int_0^t\int_0^1G(t-s,x,y)\sigma(s,y,u(s,y))W(ds,dy)\nonumber\\
&&+\int_0^t\int_0^1G(t-s,x,y)\eta(ds,dy).
\end{eqnarray}
Set
\begin{eqnarray}\label{3.22}
\bar{v}(t,x)&=&\int_0^1G(t,x,y)u(0,y)dy\nonumber\\
&&+\int_0^t\int_0^1G(t-s,x,y)f(s,y,u(s,y))dyds\nonumber\\
&&+\int_0^t\int_0^1G(t-s,x,y)\sigma(s,y,u(s,y))W(ds,dy).
\end{eqnarray}
Then $(\bar{z}(t,x):=u(t,x)-\bar{v}(t,x), \eta(dt,dx))$ solves the following random parabolic obstacle
problem:
 \begin{eqnarray}\label{3.23}\left\{
\begin{array}{l}
\frac{\partial \bar{z}(t,x)}{\partial t}-\frac{\partial^2 \bar{z}(t,x)}{\partial x^2}=\dot{\eta}(t,x), \quad x\in [0,1];\\
\bar{z}(t,x)\geq -\bar{v}(t,x);\\
\int_0^t\int_0^1(\bar{z}(s,x)+\bar{v}(s,x))\eta(ds,dx)=0,\quad t\geq 0.
\end{array}\right.
\end{eqnarray}
For very positive integer $n\geq 1$, define
$$\bar{v}^n(t)=\left(\bar{v}(t, \frac{1}{n}), ...,\bar{v}(t, \frac{n-1}{n})\right ).$$
Let $(\bar{z}^n, \bar{\eta}^n)$ be the solution of the following random  Skorohod-type  problem in $R^{n-1}$:
\begin{eqnarray}\label{3.24}\left\{
\begin{array}{l}
d\bar{z}^n(t)=n^2A^n\bar{z}^n(t)dt+d\bar{\eta}^n(t);\\
\bar{z}^n(t)\geq -\bar{v}^n(t);\\
\int_0^T<\bar{z}^n(t)+\bar{v}^n(t),d\bar{\eta}^n(t)>=0.
\end{array}\right.
\end{eqnarray}
Introduce the  continuous random field $\bar{z}^n$:
\begin{equation}\label{3.25}
\bar{z}^n(t,x):=\bar{z}_k^n(t)+(nx-k)\left ( \bar{z}_{k+1}^n(t)-\bar{z}_{k}^n(t)\right )
\end{equation}
for $x\in [\frac{k}{n}, \frac{k+1}{n})$, $k=0,...,n-1$, with $\bar{z}_0^n(t):=0$.
By Theorem 4.1, we conclude that
\begin{equation}\label{3.26}
\lim_{n\rightarrow \infty}\sup_{0\leq t\leq T, 0\leq x\leq 1}|\bar{z}^{n}(t,x)-\bar{z}(t,x)|=0
\end{equation}
almost surely.
Let
$\bar{v}^n$ denote the random field:
\begin{equation}\label{3.27}
\bar{v}^n(t,x):=\bar{v}(t, \frac{k}{n})+(nx-k)\left ( \bar{v}(t,\frac{k+1}{n})-\bar{v}(t, \frac{k}{n})\right )
\end{equation}
for $x\in [\frac{k}{n}, \frac{k+1}{n})$, $k=1,...,n-1$.
Since $\bar{v}(t,x)$ is a continuous random field with bounded moments of any order,  it is clear that for any $p\geq 1$,
\begin{equation}\label{3.28}
\lim_{n\rightarrow \infty}E[\sup_{0\leq t\leq T, 0\leq x\leq 1}|\bar{v}^{n}(t,x)-\bar{v}(t,x)|^p]=0.
\end{equation}
Set $\bar{u}^n(t,x):=\bar{v}^n(t,x)+\bar{z}^n(t,x)$. Since $u(t,x):=\bar{v}(t,x)+\bar{z}(t,x)$,  it follows from (\ref{3.26}) and (\ref{3.28}) that
\begin{equation}\label{3.29}
\lim_{n\rightarrow \infty}E[\sup_{0\leq t\leq T, 0\leq x\leq 1}|\bar{u}^{n}(t,x)-u(t,x)|^p]=0.
\end{equation}
Recall the definition of the random fields $u^n(t,x)$ defined in (\ref{3.11}) or (\ref{2.6}). To prove the theorem, i.e.,
$$
\lim_{n\rightarrow \infty}E[\sup_{0\leq t\leq T, 0\leq x\leq 1}|u^{n}(t,x)-u(t,x)|^p]=0,
$$
 in view of (\ref{3.29}) it is sufficient to show that
\begin{equation}\label{3.30}
\lim_{n\rightarrow \infty}E[\sup_{0\leq t\leq T, 0\leq x\leq 1}|\bar{u}^{n}(t,x)-u^n(t,x)|^p]=0.
\end{equation}
Applying Lemma 3.2 to the systems (\ref{3.4-1}) and (\ref{3.24}), it follows that
\begin{eqnarray}\label{3.31}
&&\sup_{0\leq t\leq T, 0\leq x\leq 1}|\bar{u}^{n}(t,x)-u^n(t,x)|\nonumber\\
&=& \sup_{0\leq t\leq T, 0\leq k\leq n-1}|\bar{u}^{n}(t,\frac{k}{n})-u^n(t,\frac{k}{n})|\nonumber\\
&=& \sup_{0\leq t\leq T, 0\leq k\leq n-1}|\bar{v}^{n}(t,\frac{k}{n})-v^{n}(t,\frac{k}{n})+(\bar{z}^n_k(t)-Z^n_k(t))|\nonumber\\
&\leq & 2\sup_{0\leq t\leq T, 0\leq k\leq n-1}|\bar{v}^{n}(t,\frac{k}{n})-v^{n}(t,\frac{k}{n})|\nonumber\\
&\leq &2\sup_{0\leq t\leq T, 0\leq x\leq 1}|\bar{v}^{n}(t,x)-v^n(t,x)|.
\end{eqnarray}
Introduce
\begin{eqnarray}\label{3.32}
\hat{v}^n(t,x)&=&\int_0^1G^n(t,x,y)u(0,k_n(y))dy\nonumber\\
&&+\int_0^t\int_0^1G^n(t-s,x,y)f(s,k_n(y),\bar{u}^n(s,k_n(y)))dyds\nonumber\\
&&+\int_0^t\int_0^1G^n(t-s,x,y)\sigma(s,k_n(y),\bar{u}^n(s,k_n(y)))W(ds,dy).
\end{eqnarray}
Recalling the expression of $v^n$ in (\ref{3.12}) we have
\begin{eqnarray}\label{3.33}
&&\hat{v}^n(t,x)-v^n(t,x)\nonumber\\
&=&\int_0^t\int_0^1G^n(t-s,x,y)
[f(s,k_n(y),\bar{u}^n(s,k_n(y)))-f(s,k_n(y),{u}^n(s,k_n(y)))]dyds\nonumber\\
&+&\int_0^t\int_0^1G^n(t-s,x,y)[\sigma(s,k_n(y),\bar{u}^n(s,k_n(y)))-\sigma(s,k_n(y),{u}^n(s,k_n(y)))]W(ds,dy).\nonumber\\
&&
\end{eqnarray}
Using the above representation, the Lipschitz continuity of the coefficients and the similar
arguments leading to the proof of (\ref{successive}) we can show that
\begin{eqnarray}\label{3.34}
E(|\hat{v}^n-v^n|^T_\infty)^p
&\leq& c(p,K,T)E\int_0^T(|\bar{u}^n-u^n|^t_\infty)^pdt.
\end{eqnarray}
Combining (\ref{3.31}) and (\ref{3.34}) we obtain that
\begin{eqnarray}\label{3.35}
&&E[(|\bar{u}^{n}-u^n|^T_\infty)^p]\nonumber\\
&\leq& cE[\sup_{0\leq t\leq T, 0\leq x\leq 1}|\bar{v}^{n}(t,x)-\hat{v}^n(t,x)|^p]\nonumber\\
&&\quad +cE[\sup_{0\leq t\leq T, 0\leq x\leq 1}|\hat{v}^{n}(t,x)-v^n(t,x)|^p]\nonumber\\
&\leq &CE[\sup_{0\leq t\leq T, 0\leq x\leq 1}|\bar{v}^{n}(t,x)-\hat{v}^n(t,x)|\nonumber\\
&&\quad +CE\int_0^T(|\bar{u}^n-u^n|^t_\infty)^pdt.
\end{eqnarray}
By the Grownwall's inequality we derive that
\begin{eqnarray}\label{3.36}
E[(|\bar{u}^{n}-u^n|^T_\infty)^p]&\leq&C(T)E[\sup_{0\leq t\leq T, 0\leq x\leq 1}|\bar{v}^{n}(t,x)-\hat{v}^n(t,x)|^p].
\end{eqnarray}
It remains to show
\begin{equation}\label{3.37}
\lim_{n\rightarrow \infty}E[\sup_{0\leq t\leq T, 0\leq x\leq 1}|\bar{v}^{n}(t,x)-\hat{v}^n(t,x)|^p]=0.
\end{equation}
From (\ref{3.27}) we deduce that
\begin{eqnarray}\label{3.38}
\bar{v}^n(t,x)&=&\int_0^1\bar{G}^n(t,x,y)u(0,y)dy\nonumber\\
&&+\int_0^t\int_0^1\bar{G}^n(t-s,x,y)f(s,y,u(s,y))dyds\nonumber\\
&&+\int_0^t\int_0^1\bar{G}^n(t-s,x,y)\sigma(s,y,u(s,y))W(ds,dy),
\end{eqnarray}
where $\bar{G}^n$ is defined as follows
\begin{equation}\label{3.39}
\bar{G}^n(t,x,y):=G(t, \frac{k}{n},y)+(nx-k)\left ( G(t,\frac{k+1}{n},y)-G(t, \frac{k}{n},y)\right )
\end{equation}
for $x\in [\frac{k}{n}, \frac{k+1}{n})$, $k=1,...,n-1$. Recall the definition of  $\varphi^n_k(x)$ in (\ref{3.1-17-1}).  It is easy to check that
$$\bar{G}^n(t,x,y)=\sum_{k=1}^{\infty}exp(-k^2\pi t)\varphi^n_k(x)\varphi_k(y),$$
and moreover, for $T>0$,
\begin{equation}\label{3.39-1}
\sup_{0\leq x\leq 1}\int_0^T\int_0^1(G(s,x,y)-\bar{G}^n(s,x,y))^2dsdy\rightarrow 0
\end{equation}
as $n\rightarrow \infty$.

\vskip 0.4cm
Now,
\begin{eqnarray}\label{3.39}
&&\hat{v}^n(t,x)-\bar{v}^n(t,x)\nonumber\\
&=&\int_0^t\int_0^1[G^n(t-s,x,y)-\bar{G}^n(t-s,x,y)]f(s,k_n(y),\bar{u}^n(s,k_n(y)))dsdy\nonumber\\
&+&\int_0^t\int_0^1[G^n(t-s,x,y)-\bar{G}^n(t-s,x,y)]\sigma(s,k_n(y),\bar{u}^n(s,k_n(y)))W(ds,dy)\nonumber\\
&+&\int_0^t\int_0^1\bar{G}^n(t-s,x,y)[f(s,k_n(y),\bar{u}^n(s,k_n(y)))-f(s,y,u(s,y))]dyds\nonumber\\
&+&\int_0^t\int_0^1\bar{G}^n(t-s,x,y)[\sigma(s,k_n(y),\bar{u}^n(s,k_n(y)))-\sigma(s,y,u(s,y))]W(ds,dy)\nonumber\\
&:=& B_1^n(t,x)+B_2^n(t,x)+B_3^n(t,x)+B_4^n(t,x).
\end{eqnarray}
We will show that each of the four terms tends to zero. In view of (\ref{3.13})) and (\ref{3.39-1}),  by the linear growth of $f$, we have
\begin{eqnarray}\label{3.40}
&&E[\sup_{0\leq t\leq T, 0\leq x\leq 1}|B_1^n(t,x)|^2]\nonumber\\
&\leq&C\left (\sup_{0\leq t\leq T, 0\leq x\leq 1}\int_0^t\int_0^1(G^n(t-s,k_n(x),y)-\bar{G}^n(t-s,k_n(x),y))^2dsdy\right )\nonumber\\
&&\quad \times \int_0^T\int_0^1(1+E[|\bar{u}^n(s,k_n(y))|^2])dsdy\nonumber\\
&\leq& C_T\sup_{0\leq t\leq T, 0\leq x\leq 1}\int_0^t\int_0^1(G^n(t-s,k_n(x),y)-\bar{G}^n(t-s,k_n(x),y))^2dsdy\rightarrow 0.
\end{eqnarray}
By the similar arguments as in the proof of  Corollary
3.4 in \cite{WA} and in the proof of Lemma 3.6 in \cite{G}, we can show that there exists a constant $K_p$ depending on  $\sup_n\sup_{0\leq t\leq T, 0\leq x\leq 1}E[|u^n(t,x)|^{2p}]$ and $\sup_{0\leq t\leq T, 0\leq x\leq 1}E[|u(t,x)|^{2p}]$ such that
\begin{equation}\label{3.41}
E[|B^n_i(t,x)-B_i^n(s,y)|^{2p}]\leq K_p(|t-s|^{\frac{1}{2}}+|x-y|)^p,
\end{equation}
for all $s, t\in [0, T], x,y\in [0, 1]$, where $i=2,3,4.$ On the other hand, for fixed
$(t,x)\in [0,T]\times [0,1]$, we have
\begin{equation}\label{3.42}
\lim_{n\rightarrow \infty}E[|B^n_i(t,x)|^{2p}]=0, \quad i=2,3,4.
\end{equation}
Let us prove (\ref{3.42}) for $B_4^n$. Other cases are similar. By Burkholder's inequality and the Lipschitz continuity of $\sigma$,
\begin{eqnarray}\label{3.43}
&&E[|B^n_4(t,x)|^{2p}]\nonumber\\
&\leq & C_pE[\left (\int_0^t\int_0^1\bar{G}^n(t-s,x,y)^2(\sigma(s,k_n(y),\bar{u}^n(s,k_n(y)))-\sigma(s,y,u(s,y)))^2
dsdy \right )^p]\nonumber\\
&\leq& C_p(\int_0^t\int_0^1\bar{G}^n(t-s,x,y)^2dsdy)^p\{\frac{1}{n}+E[\sup_{0\leq t\leq T, 0\leq x\leq 1}|\bar{u}^n(t,k_n(x))-u(t,k_n(x))|^{2p}]\nonumber\\
&& + E[\sup_{0\leq t\leq T, 0\leq x\leq 1}|u(t,k_n(x))-u(t,x)|^{2p}]\}\nonumber\\
&&\longrightarrow 0, \quad\quad \mbox{as}\quad n\rightarrow  \infty.
\end{eqnarray}
By virtue of (\ref{3.41}), (\ref{3.42}) and a standard procedure (see, e.g. \cite{Z-3}) we can deduce  that
\begin{equation}\label{3.44}
\lim_{n\rightarrow \infty}E[\sup_{0\leq t\leq T, 0\leq x\leq 1}|B^n_i(t,x)|^{2p}]=0, \quad i=2,3,4.
\end{equation}
Putting (\ref{3.39}), (\ref{3.40}) and (\ref{3.44}) together we complete the proof of (\ref{3.37}) and hence  the theorem.$\square$


\begin{thebibliography}{99}

\bibitem[DMZ]{DMZ} R. C. Dalang, C. Mueller, L. Zambotti: Hitting
properties of parabolic SPDE's with reflection. Ann. Probab.
34(4)(2006), 1423-1450.
\bibitem[DZ]{DZ} A. Debussche and L. Zambotti: Conservative stochastic Cahn-Hilliard equation with reflection. Annals Of Probability
35:5(2007) 1706-1739.
\bibitem[DP1]{MP}C. Donai-Martin, E. Pardoux: White noise driven SPDEs with reflection. Probab. Theory Relat. Fields 95, 1-24(1993).
\bibitem[DP2]{DP}C. Donai-Martin, E. Pardoux: EDPS
r$\acute{e}$fl$\acute{e}$chies et calcul de Malliavin.
(French)[SPDEs with reflection and Malliavin Calculus]. Bull. Sci.
Math. 121(5)(1997), 405-422.
\bibitem[FO]{FO}T. Funaki, S. Olla: Fluctuations for $\nabla \phi$
interface model on a wall. Stochastic Process. Appl. 94(1)(2001),
1-27.
\bibitem[G]{G}I. Gy\"o{}ngy: Lattice Approximations for Stochastic Quasi-Linear Parabolic Partial Differential Equations Driven by Space-Time White Noise I. Potential Analysis 9 (1998) 1-25.
\bibitem[G-1]{G-1}I. Gy\"o{}ngy: Lattice Approximations for Stochastic Quasi-Linear Parabolic Partial Differential Equations Driven by Space-Time White Noise II. Potential Analysis 11 (1999) 1-37.
    \bibitem[GM]{GM}I. Gy\"o{}ngy and A. Millet: Rate of convergence of space time approximations for stochastic evolution equations. Potential Analysis 30:1 (2009) 29-64.
\bibitem[LS]{LS} P.L. Lions and A.S. Sznitman: Stochastic differential equations with reflecting boundary conditions. Comm. Pure Appl. Math. XXXVII (1984) 511-537.
\bibitem[MS]{MS}
 T. Martínez and M. Sanz-Sol\`e:  A lattice scheme for stochastic partial differential equations of elliptic type in dimension $d=4$. Appl. Math. Optim. 54:3 (2006) 343–368.
\bibitem[NP]{NP}D. Nualart, E. Pardoux: White noise driven by
quaslinear SPDEs with reflection. Probab. Theory Relat. Fields
93,77-89(1992).
\bibitem[QS]{QS} L. Quer-Sardanyons and M. Sanz-Sol\`e: Space semi-discretisations for a stochastic wave equation. Potential Anal. 24:4 (2006)303–332.
    \bibitem[S]{S} L. Slomi\`n{}ski: Euler's approximations of solutions of SDEs with reflecting boundary. Stochastic Processes and Their Applications 94(2001)317-337.
\bibitem[WA]{WA} J. Walsh : An introduction to stochastic partial differential equations. In: Hennequin, P.L. (ed.), Ecole d'¨¦t¨¦ de Probabilit¨¦ de St Flour. (Lect. Notes Math., vol. 1180) Berlin Heidelberg New York: Springer
1986.
\bibitem[XZ]{XZ} T. Xu and T. Zhang : White noise driven SPDEs with reflection: existence, uniqueness and large deviation principles. Stochastic Processes and Their Applications 119:10 (2009) 3455-3470.

\bibitem[ZL-1]{ZL-1} L. Zambotti: Occupation densities for SPDEs with reflection. Annals of Probability 32:1(2004) 191-215.
\bibitem[ZL-2]{ZL-2} L. Zambotti: Integration by parts on $\delta$-Bessel bridges, $\delta>3$ and related SPDEs. Annals of Probability 31:1(2003) 323-348.
    \bibitem[ZL-3]{ZL-3} L. Zambotti: Integration by parts formulae on convex sets of paths and applications to SPDEs with reflection. Probability
    Theory and Related Fields 123:4(2002) 579-600.
\bibitem[Z-1]{Z-1} T. Zhang : White noise driven SPDEs with reflection:
strong Feller properties and Harnack inequalities. Potential Analysis 33:2 (2010) 137-151.
\bibitem[Z-2]{Z-2} T. Zhang : Large deviations for invariant measures of SPDEs with two reflecting walls. Stochastic Processes and Their Applications 122:10 (2012) 3425-3444.
\bibitem[Z-3]{Z-3} T. Zhang : Strong Convergence of Wong-Zakai Approximations of Reflected SDEs in A Multidimensional General Domain. Potential Analysis 41:3(2014) 783-815.
\end{thebibliography}
\end{document}